\documentclass[smallextended]{svjour3}       
\usepackage{latexsym}
\usepackage{amssymb}
\usepackage{amsmath}
\usepackage[mathscr]{eucal}
\usepackage{graphicx}
\usepackage{hyperref}
\renewcommand{\qed}{\hfill{\ \ \rule{2mm}{2mm}} \vspace{0.2in}}

\newcommand{\ind}{1\hspace{-2.3mm}{1}}

\begin{document}
\title{Linearity and Local Correctness in Weighted Colourings of Random Graphs}
\titlerunning{Weighted Colouring in Random Graphs}

\author{ \textbf{Ghurumuruhan Ganesan}}
\authorrunning{G. Ganesan}
\institute{IISER Bhopal,\\
\email{gganesan82@gmail.com }}
\date{}
\maketitle

\begin{abstract}
In this paper, we consider a weighted generalization of the chromatic number of  a Binomial random graph~\(G.\) We equip each edge with a random weight and then  colour the vertices in such a way that the absolute colour difference between any two adjacent vertices is at least as large as their edge weight. We show that with high probability, the weighted colouring number grows linearly with the maximum vertex degree if the edge weights have sufficiently large moments. Conversely, if the edge weight moments are unbounded then the weighted chromatic number is much larger than the maximum vertex degree, with high probability. We also obtain a sharp threshold result for locally correct weighted colourings for balanced subgraphs of~\(G.\)


\vspace{0.1in} \noindent \textbf{Key words:} Weighted colouring; random graphs; edge weight moments, locally correct colourings, sharp threshold.

\vspace{0.1in} \noindent \textbf{AMS 2020 Subject Classification:} Primary: 60C05.
\end{abstract}

\bigskip

\setcounter{equation}{0}
\renewcommand\theequation{\thesection.\arabic{equation}}
\section{Introduction}
The chromatic number~\(\chi(G)\) of a graph~\(G\) is the smallest number of colours needed for a proper colouring of~\(G\) and many upper and lower bounds for~\(\chi(G)\) exist in terms of various graph parameters like maximum vertex degree, independence number, clique number etc. (see Chapter~\(5,\) (West (2000)). For random graphs, concentration and estimates for the actual value of the chromatic number in terms of the edge probability~\(p\) has been well-studied (see for example Bollob\'as (2001)).

In this paper, we study weighted colouring in random graphs whose edge probabilities are not necessarily homogenous. Generalizing radio labelling of graphs (Chartrand et al. (2001), Korze et al. (2021))  which we call as weighted colouring number, we use the probabilistic method to estimate the weighted colouring number when the edge weights themselves are random. Such situations arise frequently in applications and we show that if the weights have sufficiently large moments, then the ``cost" incurred due to the unboundedness of weights is small.

We also consider Ramsey type properties of weighted colourings in the sense of locally correct colourings. The study of Ramsey property of random graphs is important from theoretical and application perspectives. Luczak et al. (1992) initiated the study by determining the threshold for the edge probability exponent that results in the presence of a monochromatic copy of a given graph for both vertex and edge colouring. Since then many variants and extensions have been studied. For example, R\"odl et al. (2017) extend the results of Luczak et al. (1992) and determine a double exponential upper bound for the~\(2-\)colour Folkmann numbers. Concurrently Dudek and Pralat (2017) studied monochromatic long paths in~\(3-\)colourings of random graphs and obtain optimal bounds on the length exponent. Recently Liebenau et al. (2022) established threshold properties for asymmetric Ramsey properties of~\(G\) where every colouring contains either a~\(1-\)coloured copy of~\(F\) or a~\(2-\)coloured copy of~\(H,\) where the pair~\((F,H)\) constitutes a clique and a cycle.

In this paper, we consider weighted random graphs and obtain sharp threshold for the occurrence of properly coloured copies of balanced subgraphs. We demonstrate a  ``sharp threshold" for the number of colours needed for such properly coloured copies.

The paper is organized as follows: In Section~\ref{sec_col}, we study weighted colourings and obtain estimates for the weighted chromatic number of random graphs when the edge weights have bounded moments. In Section~\ref{sec_weak}, we consider locally correct weighted colourings where we fix the number of colours and seek properly coloured subgraphs of the overall random graph.

\renewcommand\theequation{\arabic{section}.\arabic{equation}}
\section{Weighted Colourings}\label{sec_col}
Let~\(K_n\) be the complete graph on~\(n\) vertices and let~\(X_h, h \in K_n\) be independent and identically distributed (i.i.d.)  Bernoulli random variables indexed by the edge set of~\(K_n\) and satisfying
\begin{equation}\label{x_dist}
\mathbb{P}(X_f = 1)= p = 1-\mathbb{P}(X_f = 0)
\end{equation}
where~\(0 \leq p \leq 1.\) Equip the~\({n \choose 2}\) edges of~\(K_n\) with i.i.d.\ weights~\(\{w(h)\}_{h \in K_n}\) (that are also independent of~\(\{X_h\}\)), satisfying~\(w(h) \geq 1\) a.s.\ and let~\(\mathbb{P}\) denote the probability measure associated with~\(\{X(h)\}\) and~\(\{w(h)\}.\)

Let~\(G\) be the random subgraph of~\(K_n\) formed by the set of all edges satisfying~\(X_f = 1.\) For an integer~\(k \geq 1,\) a proper~\(k-\)colouring of~\(G\) is a map~\(g: V \rightarrow \{1,2,\ldots,k\}\) such that~\(g(u) \neq g(v)\) for all edges~\((u,v) \in E.\) The chromatic number~\(\chi(G)\) is the smallest integer~\(k\) such that~\(G\) admits a proper~\(k-\)colouring.

Let~\(f : V \rightarrow \{1,2,\ldots\}\) be any map from the vertex set of~\(G\) to the set of positive integers. We say that~\(f\) is a  proper weighted colouring of~\(G\) if for each edge~\(h = (u,v) \in E\) with endvertices~\(u\) and~\(v\) we have
\begin{equation}\label{rad_cond}
|f(u)-f(v)| \geq w(u,v).
\end{equation}
We define the weighted colouring number of~\(G\) as
\begin{equation}\label{grl_def}
\chi_{w} = \chi_{w}(G)  := \inf_{f} \max_{1 \leq u \leq n} f(u),
\end{equation}
where the infimum is over all proper weighted colourings of~\(G.\) Roughly speaking~\(\chi_{w}(G)\) could be thought of as a measure of the effect of weights on a proper colouring of~\(G.\) Indeed if~\(w(u,v) =1\) for all edges, then the weighted colouring number~\(\chi_{w}(G) = \chi(G),\) the chromatic number.

If~\( p \geq \frac{M\log{n}}{n}\) for some constant~\(M > 0,\) then it is well-known (see for example Theorem~\(2,\) Ganesan (2020)) that for every constant~\(\epsilon > 0\) there are constants~\(M_0,D_1,D_2 >0\) such that if~\(M > M_0,\) then
\begin{equation}\label{chi_g}
\mathbb{P}\left(\frac{D_1 np}{\log{n}} \leq \chi(G) \leq \Delta(G)+1 \leq np(1+\epsilon)\right) \geq 1-\frac{1}{n^{D_2}},
\end{equation}
where  constants do not depend on~\(n,\) throughout and~\(\Delta(G)\) is the maximum vertex degree in~\(G.\) By definition, the weighted colouring number~\(\chi_{w}(G) \geq \chi(G)\) and we have the following result regarding the linearity for~\(\chi_w(G).\)
\begin{theorem}\label{thm_crit} We have:\\
\((a)\) Suppose~\(p \geq \frac{1}{n^{\beta}}\) and~\(\mathbb{E}w^{2s+c}(h) \leq C\) for some constants\\\(0 < \beta<1, C,c>0\) and integer constant~\(s \geq 2.\) For every~\(\beta < 1-\frac{1}{s+1/2-1/2s}\) and~\(\epsilon > 0\) there is a constant~\(D > 0\) such that
\begin{equation}\label{sub_crit}
\mathbb{P}\left(\chi_w(G) \leq 2\mu np(1+\epsilon)\right) \geq 1- \frac{1}{n^{D}}
\end{equation}
where~\(\mu := \mathbb{E}w(h)\) is the expected edge weight.\\
\((b)\) Suppose~\(p = \frac{1}{n^{\beta}}\) and~\(\mathbb{P}(w(h) \geq x) \geq \frac{A}{x^{2s}}\) for some constants~\(0 < \beta < 1,A > 0\) and~\(s>1.\)
For every~\(\beta > 1-\frac{1}{2s-1},\) there are constants~\(\theta_1,\theta_2 > 0\) such that
\begin{equation}\label{super_crit}
\mathbb{P}\left(\chi_w(G) \geq n^{\theta_1} \cdot np\right) \geq 1- \exp\left(-n^{\theta_2}\right).
\end{equation}
\end{theorem}


Part~\((a)\) of Theorem~\ref{thm_crit} says that if the edge weights have bounded~\((2s)^{th}\) moments and the edge probability exponent is small enough (or equivalently if~\(G\) is sufficiently dense), then the weighted chromatic number grows linearly with the maximum vertex degree. In part~\((b),\) we obtain a partial converse of part~\((a):\) if the~\((2s)^{th}\) moment is infinite and the graph~\(G\) is sparse enough, then~\(\chi_w(G)\) is ``much" larger than the maximum degree in the sense of~(\ref{super_crit}).

We prove Theorem~\ref{thm_crit}\((a)\) using a greedy colouring strategy analogous to the maximum degree bound for unweighted chromatic number and prove part~\((b)\) by estimating the maximum weight of an edge in~\(G.\) Throughout, we use the following standard deviation estimate regarding sums of independent Bernoulli random variables. Let~\(Z_i, 1 \leq i \leq t\) be independent Bernoulli random variables satisfying~\[\mathbb{P}(Z_i = 1) = p_i = 1-\mathbb{P}(Z_i = 0).\] If~\(W_t = \sum_{i=1}^{t} Z_i\) and~\(\mu_t = \mathbb{E}W_t,\) then for any~\(0 < \epsilon < \frac{1}{2}\) we have that
\begin{equation}\label{conc_est_f}
\mathbb{P}\left(\left|W_t-\mu_t\right| \geq \epsilon \mu_t\right) \leq 2\exp\left(-\frac{\epsilon^2}{4}\mu_t\right).
\end{equation}
For a proof of~(\ref{conc_est_f}), we refer to Corollary~\(A.1.14,\) pp.~\(312,\) Alon and Spencer (2008).

\emph{Proof of Theorem~\ref{thm_crit}}: We prove~\((a)\) first and begin with a deterministic estimate that obtains a bound for the weighted chromatic number in terms of the locally averaged weights. Specifically, we show that
\begin{equation}\label{rk_est_det2}
\chi_{w}(G) \leq 1+\max_{v \in G} \sum_{u \sim v} (2w(u,v)-1)
\end{equation}
where~\(u \sim v\) denotes that~\(u\) and~\(v\) are adjacent in~\(G.\) This is a generalization of the maximum degree bound~\(\chi(G) \leq \Delta+1,\) where~\(\Delta\) is the maximum vertex degree in~\(G.\) In fact if~\(w(u,v) = 1\) for all edges, then the right side of~(\ref{rk_est_det2}) reduces to~\(\Delta+1.\) As in the proof of the maximum degree bound, we use a greedy colouring procedure to obtain~(\ref{rk_est_det2}).

Let~\(L +1:= 1+\max_{v} \sum_{u \sim v} (2w(u,v)-1),\) the right hand side of~(\ref{rk_est_det2}). We iteratively colour the vertices of~\(H\) from the set~\(\{1,2,\ldots,L+1\}.\) Pick a vertex~\(u_1\) and assign the colour~\(l(u_1) = 1.\) For the~\(i^{th}\) iteration,~\(i \geq 2,\) we pick an uncoloured vertex~\(u_i\) and let~\(l(v_1),\ldots,l(v_t)\) be the colours of its coloured neighbours. Let~\(I(v_j)\) be the set of all integers~\(l\) satisfying~\(|l-l(v_j)| \leq w(u_i,v_j)-1.\) There are~\(2w(u_i,v_j)-1\) integers in~\(I(v_j)\) and so the number of integers in~\(\bigcup_{1 \leq j \leq t} I(v_j)\) is at most~\(L.\) We therefore assign a colour~\(l(u_i) \in \{1,2,\ldots,L+1\} \setminus \bigcup_{1 \leq j \leq t}I(v_j)\) to the vertex~\(u_i.\) Continuing this process we obtain a proper weighted colouring of~\(G\) and this proves~(\ref{rk_est_det2}).

We use~(\ref{rk_est_det2}) to prove~(\ref{sub_crit}) as follows. For~\(1 \leq i \leq n\) let~\(N_i\) be the number of neighbours of the vertex~\(i\) in the graph~\(G.\) The expected number of neighbours for vertex~\(i\) is at least~\((n-1)p.\)Therefore using the concentration estimate~(\ref{conc_est_f}) with~\(\epsilon>0,\) we have that
\begin{equation}\label{ni_est}
\mathbb{P}\left(np(1-\epsilon) \leq N_i \leq  np(1+\epsilon)\right) \geq 1-e^{-Dnp},
\end{equation}
for some positive constant~\(D = D(\epsilon).\) Letting~\[E_{nei} :=\bigcap_{1 \leq i \leq n} \left\{np(1-\epsilon) \leq N_i \leq np(1+\epsilon)\right\},\] we have by the union bound that
\begin{equation}\label{nei_est}
\mathbb{P}(E_{nei}) \geq 1- ne^{-D np}.
\end{equation}

Let~\(\omega \in E_{nei}\) be a realization of~\(G\) and let~\(\mathbb{P}_{\omega}\) be the conditional probability measure given~\(\omega.\) To find an upper bound for~\(\chi_{w}(\omega),\) we use the locally averaged bound~(\ref{rk_est_det2}). For~\(1 \leq i \leq n\) let~\(J_i := \sum_{v \sim i} w(i,v) = \sum_{v=1}^{N_i}Z_v\) be the sum of the weights of edges containing~\(i\) as an endvertex, where~\(N_i\) is the degree of vertex~\(i\) in~\(\omega.\) The expected value~\(\mathbb{E}_{\omega}(J_i) = N_i\mu\) and so by the Markov inequality, we have for constant~\(\epsilon > 0\) that
\begin{eqnarray}\label{legal_blondax}
\mathbb{P}_{\omega}\left(\left|J_i-\mathbb{E}_{\omega}J_i\right| \geq \epsilon \mathbb{E}_{\omega}J_i\right) &\leq& \frac{\mathbb{E}_{\omega}\left(J_i-\mathbb{E}_{\omega}J_i\right)^{2s}}{(\epsilon\mathbb{E}_{\omega}J_i)^{2s}} \nonumber\\
&\leq& \frac{D_1}{N_i^{2s}}\mathbb{E}_{\omega}\left(J_i-\mathbb{E}_{\omega}J_i\right)^{2s}
\end{eqnarray}
for some constant~\(D_1 > 0.\) 

Expanding~\((J_i-\mathbb{E}_{\omega}J_i)^{2s} = \left(\sum_{v=1}^{N_i}(Z_v-\mathbb{E}_{\omega}Z_v)\right)^{2s},\) we see that any term of the form~\(\prod_{j}(Z_{i_j}-\mathbb{E}_{\omega}Z_{i_j})^{b_j}\) has expectation zero, unless each~\(b_j\) is even. This implies that the total number of terms with non-zero expectation is~\[\sum_{l=1}^{s}{N_i \choose l} \leq s {N_i \choose s} \leq sN_i^{s}\] by the monotonicity of the binomial coefficient for~\(N_i \geq 2s.\) Thus
\[\mathbb{E}_{\omega}\left(\sum_{v=1}^{N_i}(Z_v-\mathbb{E}_{\omega}Z_v)\right)^{2s} \leq D_2 N_i^{s}\] for some constant~\(D_2 > 0\) and plugging this into~(\ref{legal_blondax}), we get that
\[\mathbb{P}_{\omega}\left(\left|J_i-\mathbb{E}_{\omega}J_i\right| \geq \epsilon \mathbb{E}_{\omega}J_i\right) \leq \frac{D_3}{N_i^{s}}\] for some constant~\(D_3>0.\) Recalling that~\(\mathbb{E}_{\omega}(J_i) = N_i \mu\) and~\(np(1-\epsilon) \leq N_i \leq np(1+\epsilon)\) since~\(\omega \in E_{nei},\) we get
\begin{equation}\label{legal_blox}
\mathbb{P}_{\omega}\left(J_i \geq \mu(1+\epsilon)^2 np\right) \leq \frac{D_3}{N_i^{s}} \leq \frac{D_4}{(np)^{s}}
\end{equation}
for some constant~\(D_4> 0.\)



Define vertex~\(v\) to be \emph{bad} if~\(J_v > \mu np(1+\epsilon)^2\) and let~\(G_{good} \subset G\) be the graph induced by set of all vertices that are \emph{not} bad. Arguing as in the first paragraph of the proof, we obtain a proper colouring of~\(G_{good}\) using at most~\(L := 2\mu np(1+\epsilon)^2+1\) colours. Next, if~\(v_1,\ldots,v_t\) is the set of all bad vertices and~\(M_u\) is the maximum weight of an edge in~\(G\) containing~\(u\) as an endvertex, then we assign colour~\(L+2\sum_{i=1}^{j}M_{v_i}\) to vertex~\(v_j.\)  The overall resulting colouring is proper and uses at most~\(L+2\sum_{i=1}^{t}M_{v_i}\) colours
and so~\(\chi_w\) is bounded above by
\begin{eqnarray}\label{rk_est_ax}
\zeta_w &:=& 2\mu np(1+\epsilon)^2+1 + 2\sum_{i=1}^{t}M_{v_i} \nonumber\\
&=& 2\mu np(1+\epsilon)^2 +1 + 2M_{tot},
\end{eqnarray}
where~\(M_{tot} := \sum_{v} M_v \ind\left( v \text{ is bad}\right).\)

We estimate~\(M_v\ind(v \text{ is bad})\) as follows. By the moment condition in the Theorem statement and the Markov inequality, we have for~\(x \geq 1\) that
\begin{eqnarray}
  \mathbb{P}_{\omega}(M_v > x) &\leq& \mathbb{P}_{\omega}\left(\bigcup_{u \sim v}\{w(u,v) > x\}\right)  \nonumber\\
  &\leq& \sum_{u \sim v} \mathbb{P}_{\omega}(w(u,v)>x)  \nonumber\\
  &\leq&  \frac{C_1N_v}{x^{2s+c}} \nonumber\\
  &\leq& \frac{2C_1np}{x^{2s+c}} \label{thmax_ax}
\end{eqnarray}
for some constant~\( C_1>0,\) where we recall that~\(N_v\) is the number of neighbours of~\(v\) in~\(G\) and the final relation in~(\ref{thmax_ax}) is true since~\(\omega \in E_{nei}.\) Thus
\begin{eqnarray}
  \mathbb{E}M_v^{2s} &=& \int x^{2s-1}\mathbb{P}(M_v > x) dx \nonumber\\
  &\leq& 2C_1np\int_{1}^{\infty} \frac{1}{x^{1+c}} dx  \nonumber \\
   &\leq& C_2 np \nonumber
\end{eqnarray}
for some constant~\(C_2 > 0.\) By the Cauchy-Schwartz inequality and~(\ref{legal_blox}), we therefore have that
\begin{eqnarray}\label{rv_est_ax}
\mathbb{E}M_v\ind(v \text{ is bad}) &\leq& \left(\mathbb{E}M_v^{2s}\right)^{\frac{1}{2s}} \cdot \left(\mathbb{P}(v \text{ is bad })\right)^{1-\frac{1}{2s}} \nonumber\\
&\leq& (C_2np)^{\frac{1}{2s}} \cdot \left(\frac{D_5}{(np)^{s}}\right)^{1-\frac{1}{2s}} \nonumber\\
&=& \frac{C_3}{(np)^{r}}
\end{eqnarray}
where~\(r:= s - \frac{1}{2} - \frac{1}{2s}\) and~\(C_3 > 0\) is a constant.

Summing over~\(v\) in~(\ref{rv_est_ax}) and recalling that~\(M_{tot} = \sum_{v} M_v \ind(v \text{ is bad}),\) we get that
\[ \mathbb{E}M_{tot} \leq \frac{C_3n}{(np)^{r}}\] and so by Markov inequality we have for~\(\epsilon > 0\) that
\[\mathbb{P}\left(M_{tot} \geq \frac{n^{1+\epsilon}}{(np)^{r}} \right) \leq \frac{C_4}{n^{\epsilon}}\]
for some constant~\(C_4 > 0.\) Plugging this into~(\ref{rk_est_ax}), we get that
\[\mathbb{P}\left(\chi_w \leq 2\mu np(1+\epsilon)^2+1+\frac{n^{1+\epsilon}}{(np)^{r}} \right) \geq 1-\frac{C_4}{n^{\epsilon}}.\]

If we ensure that
\begin{equation}\nonumber
\frac{n^{1+\epsilon}}{(np)^{r}} \leq \epsilon np
\end{equation}
or equivalently that
\begin{equation}\label{hubba}
p \geq \left(\frac{1}{\epsilon}\right)^{\frac{1}{r+1}} \cdot \frac{1}{n^{1-\frac{\epsilon+1}{r+1}}}
\end{equation}
then we get~(\ref{sub_crit}). But by choice~\(\beta < 1-\frac{1}{r+1}\) strictly and so choosing~\(\epsilon > 0\) small enough, we see that~(\ref{hubba}) holds and this completes the proof of part~\((a)\) of the Theorem.

To prove part~\((b),\) we suppose that the event~\(E_{nei}\) defined above occurs, so that the number of edges in~\(G\) is at least~\(\frac{n^2p}{4} \geq \frac{1}{4}n^{2-\beta}\) and each edge has degree at most~\(2np = n^{1-\beta}.\) The probability that an edge~\(h\) has weight less than~\(n^{\epsilon}\) is at most~\(1-\frac{A}{n^{2s\epsilon}}\) and so the probability that the maximum edge weight is at most~\(n^{\epsilon}\) is bounded above by
\[\left(1-\frac{A}{n^{2s\epsilon}}\right)^{\frac{1}{4}n^{2-\beta}} \leq \exp\left(-\frac{A}{4}n^{2-\beta-2s\epsilon}\right) \longrightarrow 0\] if~\(\epsilon < \frac{2-\beta}{2s}\) strictly. The weight colouring number is at least as large as the maximum edge weight and so  choosing~\(1-\beta < \epsilon < \frac{2-\beta}{2s}\) (this is possible by choice of~\(\beta\)), we obtain~(\ref{super_crit}). This completes the proof of the Theorem.~\(\qed\)


\renewcommand{\theequation}{\thesection.\arabic{equation}}
\setcounter{equation}{0}
\section{Locally Correct Colourings} \label{sec_weak}
In this section, we consider a dual problem where the number of colours is fixed and we seek properly coloured subgraphs of random graphs. As before, let~\(G\) be the random graph as in~(\ref{x_dist}) and let~\(\{w(h)\}_{h \in K_n}\) be independent and identically distributed edge weights that are also independent of~\(G.\)

Let~\(\Gamma\) be a labelled deterministic graph on~\(v_0\) vertices and containing~\(e_0\) edges and let~\(V(\Gamma)\) and~\(E(\Gamma)\) respectively denote the vertex and edge sets of~\(\Gamma.\) We say that~\(\Gamma\) is \emph{balanced} if for any subgraph~\(H \subset \Gamma\) we have that~\(\frac{e(H)}{v(H)} \leq \frac{e_0}{v_0}\) where~\(e(H) = \#E(H)\) and~\(v(H) = \#V(H)\) respectively denote the number of edges and vertices in~\(H.\) Throughout~\(\#A\) denotes the cardinality of~\(A\) and we assume henceforth that~\(\Gamma\) is a balanced graph.

We say that a subgraph~\(T \subset G\) containing~\(v_0\) vertices and~\(e_0\) edges is a \emph{copy} of~\(\Gamma\) if~\(T\) is isomorphic to~\(\Gamma:\) i.e. there exists a bijection~\(g : V(T) \rightarrow V(\Gamma)\) such that~\((u,v) \in E(T)\) if and only if~\((g(u),g(v)) \in E(\Gamma).\) For constant~\(M >0,\) we say that a subgraph~\(T \subset G\) is~\(M-\)\emph{good} or simply good if
\begin{equation}\label{m_cond}
w(u,v) \leq |\theta(u)-\theta(v)| \leq M
\end{equation}
for each edge~\((u, v) \in E(T).\) Say that the colouring~\(\theta\) is \emph{good} if there is at least one good copy of~\(\Gamma\) in~\(G.\) 



Letting~\(N_{tot}\) and~\(N_{good} = N_{good}(M)\) respectively denote the total number of colourings and the number of good colourings of~\(G,\) the following result obtains a ``sharp threshold" for the occurrence of good colourings. As before, constants do not depend on~\(n.\)
\begin{theorem}\label{thm_one} Let~\(\Gamma\) be a balanced graph containing~\(v_0\) vertices and~\(e_0\) edges and suppose~\(r = n^{\theta}\text{ and } p = \frac{1}{n^{\beta}}\) for some constants~\(0 < \theta  <1\) and~\(\frac{1}{e_0} < \beta < \frac{v_0}{e_0}\) so that~\(0 < \theta_{th} := \frac{v_0-e_0\beta}{v_0-1} < 1\) strictly. There exists~\(M_0>0\) large such that for all~\(M > M_0,\) the following hold:\\
\((a)\) If~\(\theta < \theta_{th},\) then \[\mathbb{P}\left(N_{good} \geq N_{tot}\left(1-\frac{1 }{n^{C_1}}\right)\right) \geq 1- \frac{1}{n^{C_2}}\]
for some constants~\(C_1,C_2 > 0.\)\\
\((b)\) If~\(\theta > \theta_{th},\) then \[\mathbb{P}\left(N_{good} \leq \frac{N_{tot}}{n^{D_1}}\right) \geq 1- \frac{1}{n^{D_2}}\]
for some constants~\(D_1,D_2 > 0.\)
\end{theorem}
In words, we see that if~\(\theta < \theta_{th},\) then nearly all colourings are good with high probability and if~\(\theta > \theta_{th}\) then nearly all colourings are bad with high probability. Thus~\(\theta_{th}\) could be interpreted as the threshold value of~\(\theta\) for the occurrence of good colourings. 


\emph{Proof of Theorem~\ref{thm_one}}: We prove both parts in the Theorem using  method of moments and begin with some preliminary computations.  Throughout we set~\(\zeta := \theta(v_0-1) - (v_0-e_0\beta).\)

We begin  with the proof of~\((a).\) Let~\(K \geq 1\) be a large enough integer so that~\(z_K := \mathbb{P}(w(h) \leq K) > 0\) and  let~\(G_K \subset G\) be the random subgraph of~\(G\) obtained by preserving each edge~\(h\) such that~\(X_h=1\) \text{and} the weight~\(w(h) \leq K\) (see discussion following~(\ref{x_dist})). Each edge of~\(K_n\) is present in~\(G_K\) with probability~\(z_K \cdot p.\)

If~\(N_{\Gamma}(G)\) is the number of copies of~\(\Gamma\) in~\(G,\) then
\begin{equation}\label{en_gamma_est}
D_1n^{v_0}p^{e_0} \leq \mathbb{E}N_{\Gamma}(G_K)  \leq \mathbb{E}N_{\Gamma}(G) \leq D_2 n^{v_0}p^{e_0}
\end{equation}
for some constants~\(D_1,D_2 >0.\) If
\begin{eqnarray}
  E_{\Gamma} &:=& \{|N_{\Gamma}(G) - \mathbb{E}N_{\Gamma}(G)| \leq \epsilon N_{\Gamma}(G)\} \nonumber\\
   &&\;\;\; \bigcap \{|N_{\Gamma}(G_K) - \mathbb{E}N_{\Gamma}(G_K)| \leq \epsilon N_{\Gamma}(G_K)\} \label{e_gamma_def}
\end{eqnarray}
with~\(0 < \epsilon < 1\) constant, then from the Chebychev  inequality and the variance estimate for~\(N_{\Gamma}(G)\) in Lemma~\(3.5\) of Janson et al. (2000), we know that
\[\mathbb{P}\left(E^c_{\Gamma}\right)  \leq \frac{C}{\min_{H \subset \Gamma}n^{v_H}p^{e_H}}\] for some constant~\(C > 0,\) where the minimum is over all subgraphs of~\(H\) containing at least one edge and~\(v_H\) and~\(e_H\) respectively denote the number of vertices and edges in~\(H,\) respectively. Using the fact that~\(\Gamma\) is balanced, we see that~\(e_H \leq \frac{v_H e_0}{v_0}\) and so~\[\min_{H}n^{v_H}p^{e_H} \geq \min_{2 \leq i \leq v_0} n^{i}p^{ie_0/v_0} = \min_{2 \leq i \leq v_0} (n^{v_0}p^{e_0})^{i/v_0} = (n^{v_0}p^{e_0})^{2/v_0} =:n^{\nu},\] where~\(\nu > 0\) is positive since~\(\beta < \frac{v_0}{e_0}.\) Thus
\begin{equation}\label{e_gam}
\mathbb{P}\left(E^c_{\Gamma}\right)  \leq \frac{C}{n^{\nu}}.
\end{equation}

Next, let~\(X_i,  1 \leq i \leq n\) be independent and identically distributed (i.i.d.) random variables uniformly chosen from~\(\{1,2,\ldots,r\}\) that are also independent of~\(G\) and the edge weights~\(\{w(h)\}.\) Assign colour~\(X_i\) to the vertex~\(i\) and let~\(\mathbb{P}_X\) denote the distribution of~\((X_1,\ldots,X_n).\)


For a copy~\(T\) of~\(\Gamma,\) let~\(A_T\) be the event that the vertex colours~\(X_u, u \in T\) form an arithmetic progression with common difference~\(K+1.\) Also let~\(B_T\) be the event that~\(T\) is a subgraph of~\(G_K.\) If~\(A_T \cap B_T\) occurs, then~\(T\) is an~\(M-\)good copy of~\(\Gamma\) with~\(M = v_0(K+1)\) and so \[Y_{\Gamma} := \sum_{T} \ind(A_T) \ind(B_T)\] is a lower bound on the number of good copies of~\(\Gamma,\) where~\(\ind(.)\) refers to the indicator function and the summation is over all copies of~\(\Gamma\) in the \emph{complete graph}~\(K_n.\)

We now use the second moment method show that~\(Y_{\Gamma}\) is concentrated around its mean with high~\(\mathbb{P}_X-\)probability for nearly all realizations of the random subgraph~\(G.\) Indeed,
\[\mathbb{E}_XY^2_{\Gamma} = \sum_{T} \mathbb{P}_X(A_T) \ind(B_T) + \sum_{T_1 \neq T_2} \mathbb{P}_X(A_{T_1} \cap A_{T_2}) \ind\left(B_{T_1} \cap B_{T_2}\right).\]
Taking expectations with respect to~\(G\) and the edge weights, we then get that
\begin{eqnarray}\label{kalali}
\mathbb{E}\mathbb{E}_XY^2_{\Gamma} &=& \sum_{T} \mathbb{P}_X(A_T) \mathbb{P}(B_T) + \sum_{T_1 \neq T_2} \mathbb{P}_X(A_{T_1} \cap A_{T_2}) \mathbb{P}\left(B_{T_1} \cap B_{T_2}\right) \nonumber\\
&=& \sum_{T_1} \mathbb{P}_X(A_{T_1}) \mathbb{P}(B_{T_1})\left(1+ \Delta\right),
\end{eqnarray}
where
\[\Delta := \sum_{T_1 \neq T_2} \mathbb{P}_X(A_{T_2} \mid A_{T_1}) \mathbb{P}\left(B_{T_2} \mid B_{T_1}\right).\]

Let~\(T\) be any copy of~\(\Gamma\) and let~\(u\) be any vertex of~\(T.\) Conditioning on the colour~\(X_u,\) we see that~\(\mathbb{P}(A_T \mid X_u)\) lies between~\(\frac{C_1}{r^{v_0-1}}\) and~\(\frac{C_2}{r^{v_0-1}}\) for some positive constants~\(C_1,C_2\) and so taking averages we get
\begin{equation}\label{pxat_est}
\frac{C_1}{r^{v_0-1}} \leq \mathbb{P}_X(A_T)  \leq \frac{C_2}{r^{v_0-1}}.
\end{equation}
If~\(T_1\) and~\(T_2\) are distinct copies of~\(\Gamma\) that do not share any vertex, then
\begin{equation}\label{thmilp}
\mathbb{P}_X(A_{T_2} \mid A_{T_1}) = \mathbb{P}_X(A_{T_2})\text{ and }\mathbb{P}\left(B_{T_2} \mid B_{T_1}\right) = \mathbb{P}\left(B_{T_2}\right).
\end{equation}
On the other hand if~\(T_1\) and~\(T_2\) share exactly~\(1 \leq i \leq v_0-1\) common vertices, then the total number of vertices in~\(T_1 \cup T_2\) is~\(2(v_0-i)+i = 2v_0-i\) and arguing as in~(\ref{pxat_est}), we have that
\begin{equation}\label{pxat_est_ax2}
\frac{C_3}{r^{2v_0-i-1}} \leq \mathbb{P}_X(A_{T_1} \cap A_{T_2})  \leq \frac{C_4}{r^{2v_0-i-1}}
\end{equation}
for some constants~\(C_3,C_4>0.\) Combining with~(\ref{pxat_est}) we get
\begin{equation}\label{pxat_est_ax3}
\frac{C_5}{r^{v_0-i}} \leq \mathbb{P}_X(A_{T_2} \mid A_{T_1}) \leq  \frac{C_6}{r^{v_0-i}}
\end{equation}
for some constants~\(C_5,C_6 >0.\)

Next, the graph~\(\Gamma\) is balanced and so there are at most~\(\frac{i e_0}{v_0}\) edges containing both end vertices within~\(T_1 \cap T_2.\) Consequently~\[\mathbb{P}(B_{T_2} \mid B_{T_1}) \leq (pz_K)^{e_0-\frac{ie_0}{v_0}} \] where we recall that~\(z_K = \mathbb{P}(w(h) \leq K)\) is the (constant) probability that edge weight is at most~\(K.\) Since there at most~\(n^{v_0-i}\) choices for the vertices in~\(T_2 \setminus T_1,\) we get from the above discussion that
\[\Delta \leq \sum_{T_2}\mathbb{P}_X(A_{T_2})\mathbb{P}(B_{T_2})+ D_1\sum_{1 \leq i \leq v_0-1} \frac{n^{v_0-i}}{r^{v_0-i}} p^{e_0-\frac{ie_0}{v_0}}\]
for some constant~\(D_1 > 0.\) Plugging this into~(\ref{kalali}), we get that
\begin{equation}\label{kalali2}
\mathbb{E}\mathbb{E}_XY^2_{\Gamma} \leq \alpha + \alpha^2 + D_1\alpha \sum_{1 \leq i \leq v_0-1} \frac{n^{v_0-i}}{r^{v_0-i}} p^{e_0-\frac{ie_0}{v_0}},
\end{equation}
where
\begin{equation}\label{alp_def}
\alpha := \sum_{T}\mathbb{P}_X(A_{T})\mathbb{P}(B_T) = \mathbb{E}\mathbb{E}_XY_{\Gamma}.
\end{equation}

By Jensen's inequality we have that~\(\alpha^2 = (\mathbb{E}\mathbb{E}_XY_{\Gamma})^2 \leq \mathbb{E}(\mathbb{E}_XY_{\Gamma})^2\) and so from~(\ref{kalali2}) we get that
\begin{equation}\label{kalali4}
\mathbb{E} var_X(Y_{\Gamma}) \leq \alpha\left(1+D_1\sum_{1 \leq i \leq v_0-1} \frac{n^{v_0-i}}{r^{v_0-i}} p^{e_0-\frac{ie_0}{v_0}}\right).
\end{equation}
Using the lower bounds in~(\ref{pxat_est}) and~(\ref{en_gamma_est}), we get that
\begin{eqnarray}\label{alp_est}
\alpha &\geq& \sum_{T}\mathbb{P}_X(A_{T})\mathbb{P}(T \in G_K) \nonumber\\
&\geq& \frac{D_2z_K^{e_0}}{r^{v_0-1}} \mathbb{E}N_{\Gamma}(G_K) \nonumber\\
&\geq& \frac{D_3n^{v_0}p^{e_0}}{r^{v_0-1}} \nonumber\\
&=& D_3n^{|\zeta|},
\end{eqnarray}
for some constants~\(D_2,D_3 > 0,\) since~\(\theta < \frac{v_0-e_0\beta}{v_0-1}\) strictly. Similarly using the upper bounds in~(\ref{pxat_est}) and~(\ref{en_gamma_est}), we get
\begin{equation}\label{alp_est_up}
\alpha \leq \sum_{T}\mathbb{P}_X(A_{T})\mathbb{P}(T \in G)  \leq D_4n^{|\zeta|}
\end{equation}
for some constant~\(D_4 > 0.\)

Next, for~\(1 \leq i \leq v_0-1\) we have that
\[\frac{n^{v_0-i}}{r^{v_0-i}}p^{e_0-\frac{ie_0}{v_0}} = \left(\frac{n^{v_0}p^{e_0}}{r^{v_0}}\right)^{1-\frac{i}{v_0}}  \leq D_5 \left(\frac{\alpha}{r}\right)^{1-i/v_0} \leq D_5 \alpha^{1-1/v_0} = D_5\frac{\alpha}{n^{2\delta}} \] for some positive constants~\(D_5,\delta\) not depending on the choice of~\(i.\) Plugging this into~(\ref{kalali4}) we get that~\(\mathbb{E} var_X(Y_{\Gamma})  \leq \frac{D_6 \alpha^2}{n^{2\delta}}\)
for some constant~\(D_6 > 0\) and so letting~\(E_{var} := \left\{var_X(Y_{\Gamma}) \leq \frac{D_6 \alpha^2}{n^{\delta}}\right\}\) we get from the Markov inequality that
\begin{equation}\label{kalali6}
\mathbb{P}\left(E^c_{var}\right) \leq \frac{1}{n^{\delta}}.
\end{equation}

Suppose~\(E_{var} \cap E_{\Gamma}\) occurs where~\(E_{\Gamma}\) is as defined prior to~(\ref{e_gam}). For constant~\(0 < \epsilon < 1,\) we then get from the Chebychev inequality that
\begin{equation}\label{cheb_e}
\mathbb{P}_X\left(|Y_{\Gamma} - \mathbb{E}_XY_{\Gamma}| \geq \epsilon \mathbb{E}_XY_{\Gamma}\right) \leq \frac{var_X(Y_{\Gamma})}{\epsilon^2(\mathbb{E}_XY_{\Gamma})^2} \leq \frac{D_7\alpha^2}{n^{\delta}(\mathbb{E}_XY_{\Gamma})^2}
\end{equation}
for some constant~\(D_7 >0.\) From~(\ref{alp_def}) and~(\ref{pxat_est}) we have
\begin{eqnarray}\label{cheb_e3}
\mathbb{E}_XY_{\Gamma} &=&\sum_{T} \mathbb{P}_X(A_T) \ind(T \in G_K) \nonumber\\
&\geq& D_8\frac{N_{\Gamma}(G_K)}{r^{v_0-1}} \nonumber\\
&\geq& D_8\frac{(1-\epsilon)\mathbb{E}N_{\Gamma}(G_K)}{r^{v_0-1}} \nonumber\\
&\geq& D_9n^{|\zeta|}
\end{eqnarray}
for some constants~\(D_8,D_9 > 0,\) where the first relation in~(\ref{cheb_e3}) is true from the lower bound in~(\ref{pxat_est}), the second inequality in~(\ref{cheb_e3}) follows from the fact that~\(E_{\Gamma}\) occurs and the final estimate in~(\ref{cheb_e3}) is obtained by an analogous argument as in~(\ref{alp_est}). Substituting~(\ref{cheb_e3}) and the upper bound for~\(\alpha\) in~(\ref{alp_est_up}) into~(\ref{cheb_e}) we get that~\(\mathbb{P}_X\left(Y_{\Gamma} =0\right) \leq \frac{D_{0}}{n^{\delta}}\) for some constant~\(D_0 > 0.\) The estimates~(\ref{kalali6}) and~(\ref{e_gam}) then complete the proof of part~\((a)\) of the Theorem.

For part~\((b),\) we proceed in a similar manner as above with minor modifications. Let~\(J_{\Gamma}\) be the event as defined in~(\ref{e_gamma_def}) with~\(G_K\) replaced by~\(G,\) then as in~(\ref{e_gam}), we have that
\begin{equation}\label{e_gam_ax}
\mathbb{P}\left(J^c_{\Gamma}\right)  \leq \frac{C}{n^{\nu}}.
\end{equation}

As before, let~\(X_i,  1 \leq i \leq n\) be i.i.d. random variables uniformly chosen from~\(\{1,2,\ldots,r\}\) that are also independent of~\(G\) and the edge weights~\(\{w(h)\}.\) Assign colour~\(X_i\) to the vertex~\(i\) and let~\(\mathbb{P}_X\) denote the distribution of~\((X_1,\ldots,X_n).\)

For a copy~\(T\) of~\(\Gamma\) let~\(S_T\) be the event that~\(|X_u-X_v| \leq M\) for any two vertices~\(u,v \in T\) and let~\(R_T\) be the event that~\(T \in G.\) The term~\(Z_{\Gamma} = \sum_{T} \ind(S_T) \ind(R_T)\) is then an upper bound on the number of~\(M-\)good copies of~\(\Gamma.\) Let~\(T\) be isomorphic to~\(\Gamma\) and let~\(u,v\) be  vertices in~\(T.\) Given the colour~\(X_u,\) we get that~\(|X_v-X_u| \leq M\) with probability at least~\(\frac{1}{r}\) and at most~\(\frac{2M}{r}\) and since there are~\(v_0-1\) choices for~\(v,\) we get that
\begin{equation}\label{pxat_est_tt}
\frac{1}{r^{v_0-1}} \leq \mathbb{P}_X(S_T) \leq \left(\frac{2M}{r}\right)^{v_0-1}.
\end{equation}
Therefore if the event~\(J_{\Gamma}\) occurs then
\begin{eqnarray}\label{e_gam2_ax}
\mathbb{E}_X Z_{\Gamma} &\leq& N_{\Gamma}(G)\left(\frac{2M}{r}\right)^{v_0-1} \nonumber\\
&\leq& \mathbb{E}N_{\Gamma}(G)(1+\epsilon)\left(\frac{2M}{r}\right)^{v_0-1} \nonumber\\
&\leq& D_2(1+\epsilon) (2M)^{v_0-1} n^{v_0-e_0\beta - \theta(v_0-1)} \nonumber\\
&=& \frac{D_3}{n^{|\zeta|}}
\end{eqnarray}
for some constant~\(D_3 > 0,\) where the first inequality in~(\ref{e_gam2_ax}) is true by the occurrence of the event~\(J_{\Gamma}\) and the final estimate in~(\ref{e_gam2_ax}) is true  since~\(\theta > \theta_{th}\) strictly. By the Markov inequality, we therefore have that~\(\mathbb{P}_X(Z_{\Gamma} \geq 1) \leq \frac{D_4}{n^{|\zeta|}}\) for some constant~\(D_4 > 0\) and so the estimate~(\ref{e_gam_ax}) obtains part~\((b)\) of the Theorem.~\(\qed\)

\subsection*{\em Acknowledgement}
I thank Professors Rahul Roy, Thomas Mountford and C. R. Subramanian for crucial comments and also thank IMSc and IISER Bhopal for my fellowships.

\bibliographystyle{plain}

\begin{thebibliography}{10}






\bibitem{alon3} N. Alon and J. Spencer. (2008).
\newblock{\em The Probabilistic Method}.
\newblock{Wiley Interscience}.

\bibitem{boll} B. Bollob\'as. (2001).
\newblock {\em Random Graphs}.
\newblock {Cambridge University Press}.


\bibitem{chart} G. Chartrand, D. Erwin, P. Zhang, F. Harary. (2001).
\newblock{ Radio labelings of graphs}.
\newblock{\em Bull. Inst. Combin. Appl.}, \textbf{33} pp. 77-85.



\bibitem{dud} A. Dudek and P. Pralat. (2017).
\newblock{On Some Multicolor Ramsey Properties of Random Graphs}.
\newblock{\em SIAM Journal on Discrete Mathematics}, \textbf{31}, 2079--2092.


\bibitem{gan} G. Ganesan. (2020).
\newblock{Cliques and Chromatic Number in Multiregime Random
Graphs}.
\newblock{\em Sankhya},~\(doi: https://doi.org/10.1007/s13171-020-00205-4.\)



\bibitem{jan} S. Janson, T. Luczak and A. Rucinski. (2000).
\newblock{\em Random Graphs}.
\newblock{John Wiley}.


\bibitem{korz} D. Korze, Z. Shao and A. Vesel. (2021).
\newblock{New results on radio~\(k-\)labelings of distance graphs}.
\newblock{\em Discrete Applied Mathematics}, Available online, September 2021.




\bibitem{anita} A. Liebenau, L. Mattos, W. Mendonca, J. Skokan. (2022).
\newblock{Asymmetric Ramsey Properties of Random Graphs for Cliques and Cycles}.
\newblock{\em Random Structures and Algorithms}, July,  https://doi.org/10.1002/rsa.21106.

\bibitem{luc} T. Luczak, A. Ruci\'nski and B. Voigt. (1992).
\newblock{Ramsey Properties of Random Graphs}.
\newblock{\em Journal of  Combinatorial Theory Series B}, \textbf{56}, 55--68.

\bibitem{rodl} V. R\"odl, A. Ruci\'nski, M. Schacht. (2017).
\newblock{Ramsey Properties of Random Graphs and Folkmann numbers}.
\newblock{\em Discussiones Mathematicae Graph Theory}, \textbf{27}, 755--776.



\bibitem{west} D. West. (2000).
\newblock{\em Introduction to Graph Theory}.
\newblock{Prentice Hall}.

\end{thebibliography}

\end{document}